\documentclass[twoside,leqno,10pt]{amsart}
\usepackage{amsfonts}
\usepackage{amsmath}
\usepackage{amscd}
\usepackage{amssymb}
\usepackage{amsthm}
\usepackage{amsrefs}
\usepackage{latexsym}
\usepackage{bbm}
\numberwithin{equation}{section}

\begin{document}
\newtheorem{theorem}{Theorem}
\newtheorem{lemma}{Lemma}
\newtheorem*{corollary}{Corollary}
\numberwithin{equation}{section}
\newcommand{\dif}{\mathrm{d}}
\newcommand{\intz}{\mathbb{Z}}
\newcommand{\ratq}{\mathbb{Q}}
\newcommand{\natn}{\mathbb{N}}
\newcommand{\comc}{\mathbb{C}}
\newcommand{\rear}{\mathbb{R}} 
\newcommand{\prip}{\mathbb{P}}
\newcommand{\uph}{\mathbb{H}}
\newcommand{\fief}{\mathbb{F}}
\newcommand{\majorarc}{\mathfrak{M}}
\newcommand{\minorarc}{\mathfrak{m}}
\newcommand{\sings}{\mathfrak{S}}

\title{On Primes in Quadratic Progressions}
\author{Stephan Baier \and Liangyi Zhao}
\date{\today}
\maketitle

\begin{abstract}
We verify the Hardy-Littlewood conjecture on primes in
 quadratic progressions on average.  The results in the present paper
 significantly improve those of a previous paper of the authors \cite{SBLZ}.
\end{abstract}

\noindent {\bf Mathematics Subject Classification (2000)}: 11L07, 11L20, 11L40, 11N13, 11N32, 11N37 \newline

\noindent {\bf Keywords}: primes in quadratic progressions, primes represented by polynomials

\section{Introduction}
It was due to Dirichlet that any linear polynomial represents infinitely many primes provided the coefficients are co-prime.  Though long been conjectured, analogous statements are not known for any polynomial of higher degree.  G. H. Hardy and J. E. Littlewood \cite{GHHJEL} conjectured that
\begin{equation} \label{hlconj}
\sum_{n \leq x} \Lambda(n^2+k) \sim \mathfrak{S}(k) x,
\end{equation}
where $\Lambda$ is the von Mangoldt function and $\mathfrak{S}(k)$ is a
constant that depends only on $k$, as defined in
\eqref{singulardef}. Their conjecture is in an equivalent but different
form as in \eqref{hlconj}.  In fact, their conjecture is more general
than \eqref{hlconj} as it concerns the representation of primes by any quadratic polynomial that may conceivably represent infinitely many primes.  \newline

It is most note-worthy that an upper bound of the order of
magnitude predicted
by \eqref{hlconj} was proved by A. Granville and R. A. Mollin in \cite{AGRM} unconditionally uniform in
the family of quadratic polynomials, and uniform in $x$ under the
Riemann hypothesis for a certain
Dirichlet $L$-function.  Furthermore, it is shown
unconditionally in \cite{AGRM} that for large $R$ and $N$ with
$R^{\varepsilon} < N < \sqrt{R}$,
\[
 \# \left\{ n \leq N : n^2+n+A \in \prip \right\} \asymp L \left( 1, \left(
 \frac{1-4A}{\cdot} \right) \right)^{-1} \frac{N}{\log N} 
\]
holds for at least a postive proportion of integers $A$ in the range $R < A <
2R$.  They also proved in \cite{AGRM} that an asymptotic formula for
the number of prime values of $f(x)$, with $f$ belonging to certain families of quadratic
polynomials, holds for $x$ in some ranges under the assumption of the
existence of a Siegel zero for the Dirichlet $L$-function.
\newline

One may also find several results on approximations to the problem of detecting 
primes of the form $n^2+1$ in the literature. For example, Ankeny \cite{NCA} and Kubilius \cite{Kub} showed
independently that under the Riemann hypothesis for Hecke $L$-functions there
exist infinitely many primes of the form $p = m^2 + n^2$ with $n < c \log p$,
where $c$ is some positive constant. Using sieve methods, Harman and Lewis
{\cite{HaLe}} showed unconditionally that there exist infinitely many primes
of the above form with $n \le p^{0.119}$.\newline

It was established by C. Hooley \cite{CH} that if $D$ is not a perfect
square then the greatest prime factor of $n^2-D$ exceeds $n^\theta$
infinitely often if $\theta<\theta_0=1.1001\cdots$. J.-M. Deshouillers
and H. Iwaniec \cite{DesIw} improved this to the effect that $n^2+1$ has
infinitely often a prime factor greater than $n^{\theta_0-\varepsilon}$, where
$\theta_0=1.202\cdots$ satisfies $2-\theta_0-2\log(2-\theta_0)=\frac{5}{4}$.  The last-mentioned result
can also be generalized to $n^2-D$ by Hooley's arguments. \newline

Moreover, H. Iwaniec \cite{HI} showed that there are infinitely many integers
$n$ such that $n^2+1$ is the product of at most two primes.  The result
improves a previous one of P. Kuhn \cite{Kuhn} that $n^2+1$ is the product of at
most three primes for inifinitely many integers $n$ and can
extended to any irreducible polynomial $an^2+bn+c$ with $a>0$ and $c$
odd. \newline

The results mentioned in the last two paragraphs were based on sieve
methods.  It is also note-worthy that J. B. Friedlander and H. Iwaniec
\cite{FrIw2}, using results on half-dimensional sieve of H. Iwaniec
\cite{HIwan}, obtained lower bounds for the number of integers with no
small prime divisors represented by a quadratic polynomial. \newline

It is easy to see that $n^2 + 1$ represents an infinitude of
primes if and only if there are infinitely many
primes $p$ such that the fractional part of $\sqrt{p}$ is very small, namely
$< 1 / \sqrt{p}$. Balog, Harman and the first-named author \cites{3, Bal, 25} dealt with
the following related question.  Given a positive real $\lambda$ and a real
number $\theta$, for what positive numbers $\tau$ can one prove that there exist
infinitely many primes $p$ for which the inequality
$$
  \left\{ p^{\lambda} - \theta \right\} < p^{- \tau}
$$
is satisfied?  This problem in turn is related to estimating the number of primes of the form
$\left[ n^c \right]$, where $c > 1$ is fixed and $n$ runs over the
positive integers, the so-called Pyatecki\u\i-\v Sapiro
primes \cites{PiaSa}. \newline

Another approximation to the $n^2+1$ problem was given by the
authors \cite{SBLZ3}.  It was proved unconditionally that for all $\varepsilon>0$, there
exist infinitely many primes of the form $p=am^2+1$ such that $a\le
p^{5/9+\varepsilon}$, and noted that the last majorant can be taken to
be $a \leq p^{1/2+\varepsilon}$ under the
assumption of the generalized Riemann hypothesis(GRH) for Dirichlet
$L$-functions or a generalization of a conjecture of the second-named
author in \cite{LZ1} on large sieve for
square moduli.  It was also noted in \cite{SBLZ3} that under the Elliott-Halberstam \cite{PDTAEHH} conjecture
for square moduli, one can show an infinitude of primes $p=am^2+1$ with $a
\leq p^{\varepsilon}$. \newline

It is note-worthy that certain cases of the asymptotics \eqref{hlconj} would follow from the
part of another unsolved conjecture due to S. Lang and H. Trotter
{\cite{SLHT}} regarding elliptic curves.  See for
example {\cite{SLHT}} for the details.  Conjectures similar to
\eqref{hlconj} also exist for polynomials of higher degree.
Hypothesis H of A. Schinzel and W. Sierpi\'nski \cite{ASWS} gives that
if $f$ is an irreducible polynomials with integer coefficients that
is not congruent to zero modulo any prime, then $f(n)$ is prime for
infinitely many integers $n$.  P. T. Bateman and R. A. Horn
\cite{PTBRAH} gave the more explicit version, with asymptotic formula,
of the last-mentioned conjecture. \newline

We use the following standard notations and conventions in number theory throughout paper. \newline

\noindent 
$f = O(g)$ means $|f| \leq cg$ for some unspecified constant
$c>0$ which may not be the same at each occurrence. \newline
$f \ll g$ means $f=O(g)$. \newline
Following the general convention, we use $\varepsilon$ to denote a small positive constant which may not be the same at each occurrence.

\section{Statements of the Results}

The asymptotic formula in \eqref{hlconj} is studied on average by the
authors in \cite{SBLZ} and it is established that \eqref{hlconj} holds
true for almost all $k\le K$ with $x^2 (\log x)^{-A} \leq K \leq
x^2$ for any $A>0$ and noted under the assumption of GRH for Dirichlet $L$-functions that the afore-mentioned range for $K$
may be taken to be the wider range of $x^{2-\delta} \leq K \leq x^2$ for
some $\delta >0$.  In this present paper, we aim to improve the theorem
in \cite{SBLZ} and prove that \eqref{hlconj} holds for almost all
natural numbers $k \leq K$ if $x^{1+\varepsilon} \leq K \leq x^2/2$.  More precisely, we have the following in this paper.

\begin{theorem} \label{mainresult}
Suppose that $z\ge 3$. Given $B>0$, we have, for $z^{1/2+\varepsilon} \leq K \leq z/2$,
\begin{equation} \label{theoeq}
\sum_{1\le k\le  K} \left|\sum_{z<n^2+k\le 2z} \Lambda(n^2+k) -
\mathfrak{S}(k)\sum_{z<n^2+k\le 2z} 1\right|^2 \ll \frac{Kz}{(\log z)^B},
\end{equation}
where 
\begin{equation} \label{singulardef}
\mathfrak{S}(k) = \prod_{p>2} \left( 1 - \frac{\left(\frac{-k}{p}
					       \right)}{p-1} \right)
\end{equation}
with $\left( \frac{-k}{p} \right)$ being the Legendre symbol.
\end{theorem}

From Theorem~\ref{mainresult}, we deduce the following corollary.

\begin{corollary}
Given $A, B>0$ and $\mathfrak{S}(k)$ as defined in the theorem, we have, for $
z^{1/2+\varepsilon} \leq K \leq z/2$, that 
\begin{equation} \label{coroeq}
 \sum_{z<n^2+k\le 2z} \Lambda(n^2+k) = \mathfrak{S}(k) \sum_{z<n^2+k\le 2z} 1 + O \left( \frac{\sqrt{z}}{(\log z)^B} \right)
\end{equation}
holds for all natural numbers $k$ not exceeding $K$ with at most 
$O \left( K (\log z)^{-A} \right)$ exceptions.
\end{corollary}

We further note here that in \cite{SBLZ} that $k$ is set to run over
only the square-free numbers.  This unfortunate restriction is also removed in the
present paper.  Moreover, it can be easily shown, as done in section 1 of
\cite{SBLZ} that $\mathfrak{S}(k)$ converges and
\[ \mathfrak{S}(k) \gg \frac{1}{\log k} \gg \frac{1}{\log K} \gg \frac{1}{\log z}. \]
The above inequality shows that the main terms in \eqref{theoeq} and \eqref{coroeq} are indeed dominating for the $k$'s under consideration if $B>1$ and that we truly have an ``almost all'' result. \newline

Actually, we shall prove the following sharpened version of Theorem~\ref{mainresult} for short segments of quadratic progressions on average.

\begin{theorem} \label{sharperresult}
Suppose that $z\ge 3$, $z^{2/3+\varepsilon}\le \Delta\le z^{1-\varepsilon}$ and $z^{1/2+\varepsilon} \leq K \leq z/2$. Then, given $B>0$, we have
\begin{equation} \label{theoeqshort}
\int_{z}^{2z} \sum_{1\leq k\le K} \left|
\sum_{t< n^2+k \leq t+\Delta} \Lambda(n^2+k) -
\mathfrak{S}(k) \sum_{t< n^2+k \leq t+\Delta} 1 \right|^2 \dif t
\ll \frac{\Delta^2K}{(\log z)^B}.
\end{equation}
\end{theorem}  

We shall deduce Theorem~\ref{mainresult} from Theorem~\ref{sharperresult} in section 11.  Moreover, we note that under GRH, the $\Delta$-range in Theorem~\ref{sharperresult} can be extended to $z^{1/2+\varepsilon}\le \Delta\le z^{1-\varepsilon}$. This is due to the fact that under GRH, Lemmas~\ref{meanq1} and
~\ref{meansquare} hold for $\delta\ge z^{\varepsilon}$, and
Lemma~\ref{short} for $\delta\ge t^{1/2+\varepsilon}$. It is note-worthy
that for $\Delta=z^{1/2+\varepsilon}$ the segments of quadratic
progressions under consideration are extremely short; that is, they
contain only $O\left(z^{\varepsilon}\right)$ elements. \newline

Theorem~\ref{sharperresult} can be interpreted as saying that the asymptotics
$$ 
\sum_{t< n^2+k \leq t+\Delta} \Lambda(n^2+k) \sim \mathfrak{S}(k)\sum_{t< n^2+k \leq t+\Delta} 1
$$
holds for almost all $k$ and $t$ in the indicated ranges.\newline 

Unlike \cite{SBLZ}, we do not use the circle method in the present
paper.  Here our approach is a variant of the dispersion method of
J. V. Linnik \cite{Linnik}, similar to that used by H. Mikawa in the study of the twin prime problem
in \cite{Mika2}. Expanding the modulus
square in \eqref{theoeqshort}, we will get in \eqref{square} three terms $U(t)$, $V(t)$ and $W(t)$ of
which we must estimate. \newline

The cross term $V(t)$ will involve both the von Mangoldt function and
the singular series $\mathfrak{S}(k)$.  The singular series is then
split into two parts in \eqref{Vcompute}.  The first part is shown to be
small using techniques similar to those used for the analogous terms in
\cite{SBLZ}, which at the bottom invokes the large sieve for real
characters of Heath-Brown \cite{DRHB}.  The second part will give raise
to a main term using a result on primes in arithmetic progressions in
short intervals, Lemma~\ref{short}. \newline

$W(t)$, which will involve only the square of the singular series
$\mathfrak{S}(k)$, is split into two parts once again in
\eqref{singularbreak}.  One of the parts can be shown to be small using
techniques from \cite{SBLZ}.  The other will again yield a main
term using familiar estimates for character sums. \newline

The treatment of $U(t)$, the sum which will involve only the von Mangoldt
function, is the most complicated.  $U(t)$ is again decomposed into two
parts, one of which can easily be shown to be small.  Transforming the
other part of the sum, we arrive at certain congruence relations for the
summands which are detected using character sums.  The main term will,
as usual, come from the principle characters, after some
transformations.  The contribution of the non-principle characters is
once again split into two parts.  The first is disposed with the
classical large sieve and the second with second moment estimates for
certain character sums twisted with the von Mangoldt function, Lemma
\ref{meansquare}.  Further splittings are needed in the estimates of the
last-mentioned two parts to remove the dependency of certain parameters
on others. \newline

When combined, we shall discover that the three main terms mentioned
above cancel out, giving us the desired result.  The restrictions on the
sizes of $K$ and $\Delta$ are needed at various places in sections 6 -
9.  It would be highly desirable to have the results in which
$K=o(\sqrt{z})$ since in that situation the quadratic progressions under
consideration would be completely disjoint.

\section{Preliminaries}

In this section, we enumerate the lemmas needed in the proofs of the
theorems.  First, we shall use the large sieve inequality for Dirichlet characters.

\begin{lemma}[Large Sieve] \label{classls}
Let $\{ a_n \}$ be a sequence of complex numbers.  Suppose that $M \in \intz$, $N, Q \in \natn$. Then we have
$$
\sum_{Q\le q\le 2Q} \frac{1}{\varphi(q)} \sideset{}{^{\star}}\sum_{\chi \bmod q} \left| \sum_{n=M+1}^{M+N} a_n \chi (n) \right|^2 \ll \left(Q+\frac{N}{Q}\right) \sum_{n=M+1}^{M+N} |a_n|^2,
$$
where $\sideset{}{^{\star}}\sum$ henceforth denotes the sum over primitive characters to the specified modulus.
\end{lemma}

\begin{proof}
See for example \cite{HD}, \cite{PXG}, \cite{HM} or \cite{HM2} for the proof.
\end{proof}

We shall also need the following version of the large sieve for single moduli $q$.

\begin{lemma}\label{singleqls}
Let $\{ a_n \}$ be a sequence of complex numbers.  Suppose that $M \in \intz$, $N, q \in \natn$. Then we have
$$
\sideset{}{^{\star}}\sum_{\chi \bmod q} \left| \sum_{n=M+1}^{M+N} a_n \chi (n) \right|^2 \leq (q +N ) \sum_{n=M+1}^{M+N} |a_n|^2.
$$
\end{lemma}

\begin{proof}
See for example \cite{HD}, \cite{PXG}, \cite{HM} or \cite{HM2} for the proof.
\end{proof}

We shall also use the well-known estimate of Polya-Vinogradov for character sums.

\begin{lemma} [Polya-Vinogradov] \label{polyavino}
For any non-principal character $\chi \pmod{q}$ we have
\[ \left| \sum_{M <n \leq M+N} \chi (n) \right| \leq 6 \sqrt{q} \log q. \]
\end{lemma}

\begin{proof}
This is quoted from \cite{HIEK} and is Theorem 12.5 there.
\end{proof}

Furthermore, we shall use the following mean-square estimate for the von-Mangoldt function in short intervals.

\begin{lemma}\label{meanq1}
Let $z\ge 3$, $z^{1/6+\varepsilon}\le \delta \le z$ and $0<M\le\delta$. Then, for any given $C>0$, we have
$$
\int_z^{2z} \left| \sum\limits_{t<n\le t+M} \Lambda(n)\ -\ M \right|^2
 \dif t
\ll \frac{z\delta^2}{(\log z)^C}. 
$$
\end{lemma}

\begin{proof}
See Chapter 10 of \cite{HIEK} for the proof of this lemma.
\end{proof}

We shall also need the following modified version of Lemma~\ref{meanq1} for character sums with 
$\Lambda$-coefficients.

\begin{lemma}\label{meansquare}
Let $A,C>0$ be given. Suppose that $z\ge 3$, $z^{1/6+\varepsilon}\le \delta \le z$, $0<M\le \delta$ and
$2\le q\le (\log z)^A$. Then we have, for any non-principal Dirichlet character $\chi$ modulo $q$,
$$
\int_z^{2z} \left| \sum\limits_{t<n\le t+M} \Lambda(n)\chi(n) \right|^2
 \dif t \ll \frac{z\delta^2}{(\log z)^C}. 
$$
\end{lemma}

\begin{proof}
The proof goes along the same lines as Lemma~\ref{meanq1}.
\end{proof}

Furthermore, we shall use the following generalization to short intervals of the Siegel-Walfisz theorem on primes in 
arithmetic progressions.

\begin{lemma}\label{short}
Let $A,C>0$ be given. Suppose that $t\ge 3$, $t^{7/12+\varepsilon}\le \delta \le t$,
$1\le l\le (\log t)^A$ and $(a,l)=1$. Then
$$
\sum\limits_{\substack{t<n\le t+\delta\\ n\equiv a\bmod{l}}} \Lambda(n) 
= \frac{\delta}{\varphi(l)}+O\left(\frac{\delta}{(\log t)^C}\right). 
$$
\end{lemma}

\begin{proof}
For $l=1$, this is Theorem 10.5 of \cite{HIEK} with a better error term (saving of an arbitrary power of logarithm) which can be obtained from Vinogradov's widening of the classical zero-free region of the Riemann
zeta-functions. The proof for moduli $l\le (\log t)^A$ goes along the same lines by using
zero density estimates and a similar zero-free region for Dirichlet
 $L$-functions. To obtain the desired zero-free region one again uses
 Vinogradov's method together with Siegel's bound for exceptional zeros.
 See \cite{HIEK} for the details.
\end{proof}

In fact, we shall need Lemma~\ref{short} only for the range $t^{2/3+\varepsilon}\le \delta \le t$. We shall also use the following lemma on the average of $q/\varphi(4q)$.
 
\begin{lemma}\label{phi}
For $x\ge 1$ we have
$$
\sum\limits_{q\le x} \frac{q}{\varphi(4q)}=\frac{x}{2} \prod_{p>2} \left( 1+
 \frac{1}{p(p-1)} \right) + O \left(\log x\right).
$$
\end{lemma}

\begin{proof}
The proof is similar as that of (5.36) in Lemma 5.4.2. in \cite{Bru}.
\end{proof}

Finally, we shall use the following lemma on the Legendre-symbol.

\begin{lemma}\label{Legendre}
For any square-free number $l$, we have 
$$
\sum_{\substack{a\bmod{l}\\ \gcd (a,l)=1}} \ \sum\limits_{m\bmod{l}} \left( \frac{m^2-a}{l} \right) 
=\mu(l)\varphi(l).
$$
\end{lemma}

\begin{proof}
By the virtue of multiplicativity, it suffices to prove the lemma for primes $l=p$. In this case, we have
\[
\sum_{\substack{a\bmod{p}\\ \gcd (a,p)=1}} \ \sum\limits_{m\bmod{p}} \left( \frac{m^2-a}{p} \right) 
= \sum\limits_{m\bmod{p}} \ \sum_{a\bmod{p}} \left( \frac{m^2-a}{p} \right) -
\sum\limits_{m\bmod{p}} \left( \frac{m^2}{p} \right) =  -(p-1)
\]
by the orthogonality relations for Dirichlet characters. This completes the proof.
\end{proof}

\section{Preparation of the Terms}

Throughout the sequel, we assume that $z\ge 3$, $z^{1/2+\varepsilon}\le
K\le z/2$ and $z^{2/3+\varepsilon}\le \Delta\le
z^{1-\varepsilon}$. We set
\begin{equation} \label{Ldef}
L:=(\log z)^C, 
\end{equation}
where $C$ is a large positive constant. The variables $k$, $m$, $n$ denote natural numbers.\newline

We first rewrite the integrand in (\ref{theoeqshort}). 
Expanding the square, we obtain
\begin{equation} \label{square}
\sum_{1\leq k\le K} \left|
 \sum_{t< n^2+k \leq t+\Delta} \Lambda(n^2+k) -
 \mathfrak{S}(k) \sum_{t<n^2+k \leq t+\Delta} 1 \right|^2
= U(t) - 2V(t) + W(t),
\end{equation}
where
\begin{equation}\label{defU}
U(t) = \sum_{1\leq k \leq K} \mathop{\sum \sum}_{\substack{n_1, n_2 \\ t < n_1^2+k, n_2^2+k \leq t+\Delta}} \Lambda(n_1^2+k) \Lambda(n_2^2+k),
\end{equation}
\begin{equation}\label{defV}
V(t) = \sum_{1\leq k \leq K} \mathfrak{S}(k) \sum_{t < n_1^2+k \leq t+\Delta} 1  \sum_{t < n_2^2+k \leq t+\Delta} \Lambda(n_2^2+k)
\end{equation}
and
\begin{equation}\label{defW}
W(t) = \sum_{1\leq k \leq K} \mathfrak{S}^2(k) \mathop{\sum \sum}_{\substack{n_1, n_2 \\ t < n_1^2+k, n_2^2+k \leq t+\Delta}} 1.
\end{equation}

As mentioned in the introduction, we shall develop asymptotic formulas
for $U(t)$, $V(t)$ and $W(t)$, and the main terms will cancel out, giving
the desired result.

\section{Decomposition of $U(t)$}
We aim to derive an upper bound of correct order of magnitude for 
$$
\int\limits_{z}^{2z} U(t) \dif t.
$$  
We note that the average order of $U(t)$ may be expected to be 
$$
\sim \frac{\Delta^2K}{z}.
$$
Therefore, on average, the quantity
$$
E:=\frac{\Delta^2K}{z\log^{B}z}, \; \mbox{with} \; B>0,
$$
should be small compared to $U(t)$. In the sequel, we will make frequent use of the quantity $E$ to bound error terms.\newline

We now decompose $U(t)$ into two parts
\begin{equation}\label{deco}
U(t)=\sum_{1\leq k \leq K} \mathop{\sum \sum}_{\substack{n_1, n_2 \\ t < n_1^2+k, n_2^2+k \leq t+\Delta\\ 
|n_1^2-n_2^2|\le \Delta/L}} \Lambda(n_1^2+k) \Lambda(n_2^2+k) 
+ 2\sum_{1\leq k \leq K} \mathop{\sum \sum}_{\substack{n_1, n_2 \\ t < n_1^2+k, n_2^2+k \leq t+\Delta\\ 
n_1^2-n_2^2> \Delta/L}} \Lambda(n_1^2+k) \Lambda(n_2^2+k),
\end{equation}
where $L$ is defined in \eqref{Ldef}. The first sum on
the right-hand side of \eqref{deco} is easily seen to be $O(E)$. 
Re-writing $m_i=n_i^2+k$ with $i=1$ and $2$, the second sum on the right-hand side of 
\eqref{deco} is
$$
 \tilde{U}(t)= 2\mathop{\sum \sum}_{\substack{t < m_1, m_2 \leq
 t+\Delta\\ m_1-m_2>\Delta/L}} \Lambda(m_1) \Lambda(m_2)
 \sum_{\substack{n_1,n_2\\ m_1-m_2=n_1^2-n_2^2 \\ 1 \leq m_1-n_1^2 \leq K}} 1.
$$
It suffices to consider only the case when both $m_1$ and $m_2$ are odd at the
cost of a small error of size $\ll E$.
Now we set
\[
 q=(n_1-n_2)/2 \ \ \ \ \ \ \mbox{and} \ \ \ \ \ \ r=(n_1+n_2)/2.
\]
In the case that $m_1$ and $m_2$ are both odd, $m_1-m_2$ is even and
hence the condition
\[
  m_1-m_2=n_1^2-n_2^2
\]
implies that $n_1$ and $n_2$ are of the same parity. Therefore  $q$ and $r$ are integers in this case. Moreover, the condition $m_1-m_2=n_1^2-n_2^2$ is equivalent to 
\[
  m_1-m_2=4qr.
\]
Now $\tilde{U}(t)$ becomes
$$
\tilde{U}(t) =  2 \mathop{\sum \sum}_{\substack{t < m_1, m_2 \leq
 t+\Delta \\ m_1-m_2 > \Delta/L}} \Lambda(m_1) \Lambda(m_2) \sum_{\substack{q,r\in \mathbbm{N}\\ m_1-m_2=4qr \\ 0
 < m_1-(q+r)^2 \leq K \\ |q|<|r| }} 1 \ +\ O(E).
$$
We note that
$$
\sqrt{z}\le \sqrt{m_1-K}+\sqrt{m_2-K}\le 2r=n_1+n_2\le \sqrt{m_1}+\sqrt{m_2}\le 4\sqrt{z}
$$
if $z$ is 
sufficiently large. Therefore, the variable $q$ satisfies the condition 
$$
D_1\leq q=\frac{m_1-m_2}{4r}\le D_2,
$$
where
$$
D_1:=\frac{\Delta}{8L\sqrt{z}},\ \ \ \ \ \ 
D_2:=\frac{\Delta}{2\sqrt{z}}.
$$ 
Moreover, if $m_1>m_2$, the condition
\[
 0 < m_1-\left( q+ r \right)^2 = m_1 - \left( q+ \frac{m_1-m_2}{4q} \right)^2 \leq K
\]
holds if and only if
\[
 m_1 - 4q(\sqrt{m_1}-q) < m_2 \leq m_1 - 4q(\sqrt{m_1-K}-q).
\]
Now we set
\[
 \mathcal{I}(t,m,q) = \left( m - 4q(\sqrt{m}-q) , 
m-4q(\sqrt{m-K}-q) \right]\cap (t,t+\Delta] .
\]
Then $\tilde{U}(t)$ is majorized by
\begin{equation} \label{nomin}
 2 \sum_{D_1\leq q \leq D_2}\ \sum_{t<m_1\leq t+\Delta}
 \Lambda(m_1) \sum_{\substack{ 
 m_2 \in \mathcal{I}(t,m_1,q)\\ m_2\equiv m_1 \bmod{4q}}} \Lambda(m_2)\ + \ O(E).
\end{equation}

Due to the presence of $\Lambda$, the contribution of $m_1$ and $m_2$ in \eqref{nomin} that are not prime to
$4q$ is small, $O(E)$ with an absolute implied constant.  For the $m_1$ and $m_2$ that are prime to $4q$, we 
use Dirichlet characters to detect the congruence relation in
\eqref{nomin}, and this part becomes
\begin{equation} \label{insertchar}
2 \sum_{D_1\le q\le D_2} \frac{1}{\varphi(4q)} \sum_{\chi \bmod{4q}}\ \sum_{t<m_1\leq t+\Delta} \Lambda(m_1) \chi(m_1) \sum_{m_2 \in \mathcal{I}(t,m_1,q)} \Lambda(m_2) \overline{\chi} (m_2).
\end{equation}
\newline
The main term in \eqref{insertchar} comes from the principal characters. Up to a small error of size $\ll E$, this main term amounts to
\begin{equation} \label{main}
M(t)=2 \sum_{D_1\le q\le D_2} \frac{1}{\varphi(4q)} \sum_{t<m_1\leq t+\Delta}
 \Lambda(m_1) \sum_{m_2 \in \mathcal{I}(t,m_1,q)} \Lambda(m_2).
\end{equation}
We will deal with $M(t)$ later in section 8. Up to a small error of size $\ll E$, the remaining part of \eqref{insertchar} can be rewritten in the form
$$
2 \sum_{d\leq 2D_2}\ \sum_{\substack{\max\{2,4D_1/d\}\leq q_1 \leq 4D_2/d\\
4|q_1d}} \ \frac{1}{\varphi(q_1d)}\ \sideset{}{^{\star}}\sum_{\chi \bmod{q_1}} \
\sum_{t<m_1\leq t+\Delta} \Lambda(m_1) \chi(m_1) \sum_{m_2\in \mathcal{I}(t,m_1,q_1d/4)} \Lambda(m_2) \overline{\chi} (m_2)=:
2F(t),
$$
say.  We write
$$
F(t)=\sum_{d\leq D_1/L} \cdots \ \ + \sum_{D_1/L<d\leq 2D_2} \cdots \ \
=F_1(t)+F_2(t), \; \mbox{say}.
$$
We note that the expression $F_2(t)$ involves only small moduli $q_1\ll L^2$, whereas the moduli $q_1$ contained in $F_1(t)$ satisfy the inequality $q_1\ge 4L$.\newline

\section{Estimation of $F_1(t)$}

In this section, we shall show that $F_1(t)$ is an error term, {\it i.e.} $F_1(t)\ll E$. 
To separate the sums over $m_1$ and $m_2$ contained in $F_1(t)$, we split the ranges of summation for $q_1$ and $m_1$ into certain subintervals and then approximate the range $\mathcal{I}(t,m_1,q_1d/4)$ of summation for $m_2$ suitably. More in particular, we split the summation interval $4D_1/d\leq q_1 \leq 4D_2/d$ into $O(\log z)$ dyadic intervals $(Q,2Q]$ and then 
split the summation interval $t<m_1\le t+\Delta$ into 
$O(\Delta L/T)$ subintervals $(s,s+M]$ of length $M\le T/L$, where
$$
T := \frac{QdK}{\sqrt{z}}.
$$
We note that the inequality $T/L\le \Delta$ is always satisfied since $Qd\le \Delta/\sqrt{z}$ and $K\le z$. Now, when $Q< q_1\le 2Q$ and $s<m_1\leq s+M$, we  
replace $\mathcal{I}(t,m_1,q_1d/4)$ with $\mathcal{I}(t,s,q_1d/4)$ in the range of
summation of $m_2$. The error $R(t)$ caused by this change turns out to be small, {\it i.e.} it is $\ll E$. In the following, we indicate how the latter can be proved, but we skip the details.\newline

The error term $R(t)$ in question is a sum over $d$, $q_1$, the primitive characters $\chi$ modulo $q_1$, 
$m_1$, and $m_2$. Here the inner-most sum over $m_2$ ranges over small intervals of length $\ll T/L$ since 
from $|m_1-s| \le T/L$, 
$z^{1/2+\varepsilon}\le K\le z/2$, $\Delta\le z^{1-\varepsilon}$ and 
$Qd\ll \Delta/\sqrt{z}$ it follows that
$$
|\{m_1 - 4q_1d(\sqrt{m_1}-q_1d)\}-\{s - 4q_1d(\sqrt{s}-q_1d)
\}|\ll T/L
$$
and 
$$
|\{m_1 - 4q_1d(\sqrt{m_1-K}-q_1d)\}-\{s - 4q_1d(\sqrt{s-K}-q_1d)\}|\ll T/L.
$$ 
We note that, in contrast, the length of the interval 
$$\left( m_1 - 4q_1d(\sqrt{m_1}-q_1d),
m_1 - 4q_1d(\sqrt{m_1-K}-q_1d) \right]$$ 
is $\gg T$. Now we estimate the sums in $R(t)$ trivially.  After a short
computation, we arrive at the desired bound $R(t)\ll E$. \newline

The remaining task in this section is to establish an estimate for 
\begin{equation} \label{expression}
\sum_{Q< q_1 \leq 2Q} \frac{1}{\varphi(q_1d)} 
\left| \ \sideset{}{^{\star}} \sum_{\chi \bmod{q_1}} 
\sum_{s<m_1\leq s+M} \Lambda(m_1) \chi(m_1) \sum_{m_2\in \mathcal{I}(t,s,q_1d/4)} \Lambda(m_2) \overline{\chi} (m_2)\right|.
\end{equation}
If we can show that \eqref{expression} satisfies the non-trivial bound
\begin{equation} \label{aim}
\ll \frac{T^2\log^2 z}{\varphi(d)L^{3/2}},
\end{equation}
then it can now be easily deduced that
$$
F_1(t)\ll E,
$$
as desired. Using Cauchy's inequality and the inequality 
$$
\frac{1}{\varphi(q_1d)}\le \frac{1}{\varphi(q_1)}\cdot \frac{1}{\varphi(d)},
$$
\eqref{expression} is bounded by
\begin{equation} \label{CS}
\ll \frac{(S_1S_2)^{1/2}}{\varphi(d)},
\end{equation}
where 
$$
S_1=\sum_{Q< q_1 \leq 2Q} \frac{1}{\varphi(q_1)} 
\sideset{}{^{\star}}\sum_{\chi \bmod{q_1}} \left| 
\sum_{s<m_1\leq s+M} \Lambda(m_1) \chi(m_1)\right|^2
$$
and 
$$
S_2=\sum_{Q< q_1 \leq 2Q} \frac{1}{\varphi(q_1)} 
\sideset{}{^{\star}}\sum_{\chi \bmod{q_1}} \left|
\sum_{m_2\in \mathcal{I}(t,s,q_1d/4)} \Lambda(m_2) \overline{\chi} (m_2)\right|^2.
$$
Using the large sieve inequality, Lemma~\ref{classls}, we obtain 
\begin{equation} \label{S1}
S_1\ll \left(Q+\frac{T}{QL}\right)\frac{T\log^2 z}{L}.
\end{equation}
From Lemma~\ref{singleqls}, we deduce that
\begin{equation} \label{S2}
S_2\ll \left(Q+T\right)T\log^2 z.
\end{equation}
We note that
\begin{equation} \label{Qcond}
L\le Q\le \frac{T}{L^2}.
\end{equation}
The first inequality in \eqref{Qcond} comes from the fact that the moduli $q_1$ in the expression $F_1(t)$
satisfy the inequality $q_1\ge L$. The second inequality in \eqref{Qcond} follows from
the definition of $T$ and $K\ge z^{1/2+\varepsilon}$. From \eqref{CS},
\eqref{S1}, \eqref{S2} and \eqref{Qcond}, we obtain that
\eqref{expression} satisfies the majorant in \eqref{aim}, as desired. \newline

\section{Treatment of $F_2(t)$}
Next, we turn to the term $F_2(t)$.
We recall that the expression $F_2(t)$ involves only small moduli $q_1$. More in particular, we have that $q_1\ll L^2$. Fix any $q_1$ satisfying this inequality. 
Similarly as in the estimation of $F_1(t)$, we split the summation interval $t<m_1\le t+\Delta$ into 
$O(\Delta L/T)$ subintervals $(s,s+M]$ of length $M\le T/L$, where now we set
$$
T:= \frac{q_1dK}{\sqrt{z}}.
$$
As before, we  
replace $\mathcal{I}(t,m_1,q_1d/4)$ with $\mathcal{I}(t,s,q_1d/4)$ in the range of
summation of $m_2$ at the cost of a small error whose total contribution to $F_2(t)$ is $\ll E$. Moreover, we put
$$
s=t+\sigma.
$$
We aim to show that $F_2(t)$ is small on average, {\it i.e.} 
\begin{equation} \label{F2}
\int_z^{2z} F_2(t) \dif t \ll zE=\frac{\Delta^2K}{\log^{B}z}.
\end{equation}
To establish \eqref{F2}, it now suffices to prove a non-trivial estimate of the form
\begin{equation} \label{thegoal}
\int_z^{2z} \left|
\sum_{t+\sigma<m_1\leq t+\sigma+M} \Lambda(m_1) \chi(m_1) \sum_{m_2\in \mathcal{I}(t,t+\sigma,q_1d/4)} \Lambda(m_2) \overline{\chi} (m_2)\right|\dif t
\ll \frac{zT^2\log z}{L^2}
\end{equation}
for any fixed $\sigma$ with $0\le \sigma\le \Delta$ and any primitive character $\chi$ with
conductor $q_1$. Using Cauchy's inequality, the left-hand side of \eqref{thegoal} is bounded by
\begin{equation} \label{CI}
\ll (I_1I_2)^{1/2},
\end{equation}
where 
$$
I_1=\int_z^{2z} \left| \sum_{t+\sigma<m_1\leq t+\sigma+M} \Lambda(m_1) \chi(m_1) 
\right|^2\dif t
$$
and 
$$
I_2=\int_z^{2z} \left| \sum_{m_2\in \mathcal{I}(t,t+\sigma,q_1d/4)} \Lambda(m_2) \overline{\chi} (m_2)\right|^2\dif t.
$$
Taking into account that 
$$
\frac{T}{L}=\frac{q_1dK}{\sqrt{z}L}\ge \frac{D_1K}{\sqrt{z}L} \gg 
\frac{\Delta K}{zL^2}\gg z^{1/6+\varepsilon},
$$
we deduce from Lemma~\ref{meansquare} that
\begin{equation} \label{I1}
I_1\ll \frac{zT^2}{L^4\log^2 z}.
\end{equation}
Now it already suffices to estimate $I_2$ trivially by 
\begin{equation} \label{I2}
I_2\ll  zT^2\log^2 z.
\end{equation}
Combining \eqref{CI}, \eqref{I1} and \eqref{I2}, we obtain \eqref{thegoal}.

\section{Contribution of the main term $M(t)$}

Now we want to derive an asymptotic estimate for 
$$
\int_z^{2z} M(t)\dif t,
$$
the integral of the main term $M(t)$, defined in \eqref{main}.  
In a similar way as we established \eqref{F2} in the previous section,
it can be shown that
\begin{equation} \label{app}
\int_z^{2z} \left( M(t)-\tilde{M}(t) \right) \dif t \ll zE=\frac{\Delta^2K}{\log^{B}z}
\end{equation}
with 
$$
\tilde{M}(t)=2 \sum_{D_1\le q\le D_2} \frac{1}{\varphi(4q)} \sum_{t<m_1\leq t+\Delta}
\sum_{m_2 \in \mathcal{I}(t,m_1,q)} 1,
$$
where here we use Lemma~\ref{meanq1} instead of Lemma~\ref{meansquare}. It thus suffices to
establish an asymptotic estimate for $\tilde{M}(t)$. We shall show that 
\begin{equation} \label{finalM}
\tilde{M}(t)= \frac{\Delta^2 K}{4t} \prod_{p>2} \left( 1+
 \frac{1}{p(p-1)} \right) + O \left(E\right).
\end{equation}\newline

First, we change the order of summation, thus obtaining
$$
\tilde{M}(t)=2 \sum_{t<m_1\leq t+\Delta}\ \sum_{D_1\le q\le D_2} \frac{1}{\varphi(4q)} 
\sum_{m_2 \in \mathcal{I}(t,m_1,q)} 1.
$$
Now we make two replacements, each at the cost of an error of size $\ll E$. First we replace the summation interval for $q$  by $1\le q\le (m_1-t)/(4\sqrt{m_1})$, and second we replace the summation interval $\mathcal{I}(t,m_1,q)$ for $m_2$ by
$$ 
\left(m_1-4q(\sqrt{m_1}-q),m_1-4q(\sqrt{m_1}-q)+2qK/\sqrt{m_1}\right].
$$ 
Thus, we obtain
\begin{eqnarray*}
\tilde{M}(t) &=& 2 \sum_{t<m_1\leq t+\Delta}\ \sum_{q\le (m_1-t)/(4\sqrt{m_1})} \frac{1}{\varphi(4q)} \cdot \frac{2qK}{\sqrt{m_1}} + O(E)\\
&=& 4K \sum_{t<m_1\leq t+\Delta} \frac{1}{\sqrt{m_1}} \sum_{q\le (m_1-t)/(4\sqrt{m_1})} \frac{q}{\varphi(4q)} + O(E).
\end{eqnarray*}
Using Lemma~\ref{phi}, the above is
\begin{eqnarray*}
&=& \frac{K}{2} \prod_{p>2} \left( 1+
 \frac{1}{p(p-1)} \right) \sum_{t<m_1\leq t+\Delta} \frac{m_1-t}{m_1} + O(E)\\
&=& \frac{\Delta^2K}{4t} \prod_{p>2} \left( 1+
 \frac{1}{p(p-1)} \right) + O(E),
\end{eqnarray*}
which completes the proof of \eqref{finalM}.

\section{Contribution of $V(t)$}

We now consider $V(t)$. Expanding the Euler product that defines $\mathfrak{S}(k)$, and approximating the sum $$\sum_{t < n_1^2+k \leq t+\Delta} 1$$ by $\Delta/(2\sqrt{t})$, we 
obtain
$$
V(t) = \frac{\Delta}{2\sqrt{t}} \sum_{1\leq k \leq K}\ \sum_{\substack{l=1 \\ 2 \nmid l}}^{\infty} \frac{\mu(l)}{\varphi(l)} \left( \frac{-k}{l} \right) \sum_{t < n^2+k \leq t+\Delta} \Lambda (n^2+k)+O(E).
$$
The right-hand side of the above can be re-written as
\begin{equation} \label{Vcompute}
\begin{split}
 \frac{\Delta}{2\sqrt{t}} \sum_{1\leq k \leq K}\ & \sum_{\substack{l > L \\ 2 \nmid l}} \frac{\mu(l)}{\varphi(l)} \left( \frac{-k}{l} \right) \sum_{t < n^2+k \leq t+\Delta} \Lambda (n^2+k) \\
& +\frac{\Delta}{2\sqrt{t}} \sum_{\substack{l \leq L \\ 2\nmid l}} \frac{\mu(l)}{\varphi(l)} \sum_{1 \leq k \leq K}\ \sum_{t< n^2+k \leq t+\Delta} \Lambda(n^2+k) \left( \frac{-k}{l} \right) \  +\  O(E).
\end{split}
\end{equation}
The first term in \eqref{Vcompute} is, by Cauchy's inequality,
\begin{eqnarray*}
& \ll & \frac{\Delta}{\sqrt{z}} \left( \sum_k \left| \sum_{n} \Lambda(n^2+k) \right|^2 \right)^{1/2} \times \left( \sum_k \left| \sum_{\substack{l > L \\ 2 \nmid l}} \frac{\mu(l)}{\varphi(l)} \left( \frac{-k}{l} \right) \right|^2 \right)^{1/2} \\ 
& \ll & \frac{\Delta}{\sqrt{z}} \cdot \frac{\Delta \sqrt{K}\log z}{\sqrt{z}}\cdot  \frac{\sqrt{K}}{(\log z)^{B+1}}=E,
\end{eqnarray*}
where we have used the same techniques as in section 5 of \cite{SBLZ}, with the
following observation
\[
 \sum_k \left| \sum_{\substack{l > L \\ 2 \nmid l}}
 \frac{\mu(l)}{\varphi(l)} \left( \frac{-k}{l} \right) \right|^2 =
 \sum_{s \leq \sqrt{K}}\ \sum_{\substack{k\leq K/s^2 \\ \mu(k) \neq 0}}
 \left| \sum_{\substack{l > L \\ 2 \nmid l}}
 \frac{\mu(l)}{\varphi(l)} \left( \frac{-k}{l} \right) \right|^2.
\]

We now deal with the inner double sum in the second term of \eqref{Vcompute}.  It is
\begin{eqnarray}
\nonumber & = & \sum\limits_{t<b\le t+\Delta} \Lambda(b) \sum_{\substack{n\\ 1 \leq b-n^2 \leq t+\Delta}} \left( \frac{n^2-b}{l} \right)\\ \label{new}
& = & \sum_{a\bmod{l}} \ \sum\limits_{m\bmod{l}} \left( \frac{m^2-a}{l} \right) \sum\limits_{\substack{t<b\le t+\Delta\\ b\equiv a\bmod{l}}} \Lambda(b) \sum_{\substack{n\equiv m\bmod{l}\\ 1 \leq b-n^2 \leq K}} 1.
\end{eqnarray}
Approximating the inner-most sum of \eqref{new} by $K/(2l\sqrt{t})$, and estimating the sum of $\Lambda(b)$ by Lemma~\ref{short} upon noting that $\Delta \ge z^{2/3+\varepsilon}$, we transform \eqref{new} into
$$
 = \frac{\Delta K}{2\varphi(l)l\sqrt{t}}\ \sum_{\substack{a\bmod{l}\\ (a,l)=1}} \ \sum\limits_{m\bmod{l}} \left( \frac{m^2-a}{l} \right) + O\left(\frac{\Delta K}{\sqrt{z}(\log z)^{B+1}}\right).
$$
By Lemma~\ref{Legendre}, the above is 
$$
= \frac{\Delta K}{2\sqrt{t}} \cdot \frac{\mu(l)}{l} + O\left(\frac{\Delta K}{\sqrt{z}(\log z)^{B+1}}\right)
$$
if $l$ is square-free which can be assumed due to the presence of $\mu(l)$ in \eqref{Vcompute}.
Now, combining everything, we obtain
\begin{equation} \label{finalV}
V(t)= \frac{\Delta^2 K}{4t} \sum_{\substack{l \leq L \\ 2 \nmid l}}
 \frac{\mu(l)^2}{l \varphi(l)} + O(E) =  
 \frac{\Delta^2 K}{4t} \prod_{p>2} \left( 1+
 \frac{1}{p(p-1)} \right) + O(E).
\end{equation}

\section{Contribution of $W(t)$}

We now consider $W(t)$.  We may approximate $W(t)$ by
\begin{equation} \label{Wcompute}
  W(t) = \frac{\Delta^2}{4t} \sum_{1\leq k \leq K} \mathfrak{S}^2 (k)  + O \left(E\right).
\end{equation}
Also we have that
\begin{equation} \label{singularbreak}
  \mathfrak{S} (k) = \sum_{\substack{1 \leq l \leq L \\ 2 \nmid l}}
  \frac{\mu(l)}{\varphi(l)} \left( \frac{-k}{l} \right) + \sum_{\substack{l > L \\ 2 \nmid l}}
  \frac{\mu(l)}{\varphi(l)} \left( \frac{-k}{l} \right).
\end{equation}
First we have
\[
\sum_k \sum_{\substack{l > L \\ 2 \nmid l}}
  \frac{\mu(l)}{\varphi(l)} \left( \frac{-k}{l} \right) \ll \frac{K}{(\log z)^{B}},
\]
again using the same techniques as in section 5 \cite{SBLZ}.  Expanding the sum over $k$ of the square of the first term in \eqref{singularbreak}, we get
\begin{equation} \label{so}
\sum_{1\leq k \leq K} \left| \sum_{\substack{1 \leq l \leq L \\ 2 \nmid l}}
  \frac{\mu(l)}{\varphi(l)} \left( \frac{-k}{l} \right) \right|^2 =
  \sum_{1\leq l \leq L} \frac{\mu^2(l)}{\varphi^2(l)} \left( K \frac{\varphi(l)}{l} + O(l) \right) + \mathop{\sum \sum}_{\substack{1 \leq l_1, l_2 \leq L \\ l_1 \neq l_2}} \frac{\mu(l_1) \mu(l_2)}{\varphi(l_1) \varphi(l_2)} \sum_{1\le k\le K} \left( \frac{-k}{l_1l_2} \right).
\end{equation}
Using the inequality of Polya-Vinogradov, Lemma~\ref{polyavino}, we deduce that the
second term on the right-hand side of \eqref{so} is bounded by 
\[
 \ll \frac{K}{(\log z)^{B}}.
\]
Completing the sum 
$$\sum_{1\leq l \leq L} \frac{\mu^2(l)}{l\varphi(l)}$$
contained in the first term on the right-hand side of \eqref{so}, and combining everything,
we arrive at the estimate
\begin{equation} \label{finalW}
W(t)=  
 \frac{\Delta^2 K}{4t} \prod_{p>2} \left( 1+
 \frac{1}{p(p-1)} \right) + O \left(E\right).
\end{equation}

\section{Proofs of Theorems~\ref{mainresult} and \ref{sharperresult}}

We first prove Theorem~\ref{sharperresult}. 
\begin{proof} (of Theorem~\ref{sharperresult})
Combining the estimates derived in sections 5-8, we obtain the following bound
\begin{equation} \label{finalU}
\int_z^{2z} U(t)\dif t \le \Delta^2 K \cdot \frac{\log 2}{4} \prod_{p>2} \left( 1+
 \frac{1}{p(p-1)} \right) + O\left(\frac{\Delta^2K}{\log^{B}z}\right).
\end{equation}
Using \eqref{finalU}, \eqref{finalV} and \eqref{finalW}, and taking into account that 
$U(t)-2V(t)+W(t)$ is positive by \eqref{square}, we deduce that
$$
\int_z^{2z} (U(t)-2V(t)+W(t)) \dif t\ll \frac{\Delta^2K}{(\log z)^{B}}.
$$
From the above and \eqref{square}, we obtain the desired estimate
 \eqref{theoeqshort}.
\end{proof}

Now we derive Theorem~\ref{mainresult} from
Theorem~\ref{sharperresult}. 
\begin{proof} (of Theorem~\ref{mainresult})
We assume that $x^{2/3+\varepsilon}\le \Delta\le x^{1-\varepsilon}$ and write
\begin{eqnarray} \label{start}
& & \sum_{z<n^2+k \le 2z} \Lambda(n^2+k) -
\mathfrak{S}(k)\sqrt{z} \\ &=& \sum_{z<n^2+k \le 2z} \Lambda(n^2+k) -
\mathfrak{S}(k)\sum_{z<n^2+k \leq 2z} 1 + O\left(\frac{\sqrt{z}}{(\log z)^B}\right) \nonumber\\ 
&=&
\frac{1}{\Delta} \int_{z}^{2z} \left(
\sum_{t< n^2+k \leq t+\Delta} \Lambda(n^2+k) -
\mathfrak{S}(k) \sum_{t< n^2+k \leq t+\Delta} 1 \right) \dif t
+ O\left(\frac{\sqrt{z}}{(\log z)^B}\right),\nonumber
\end{eqnarray}
where we have used that $\mathfrak{S}(k)\ll \log 2k$ as was shown in
 section 1 of \cite{SBLZ}.  Using \eqref{start} and Cauchy's inequality, the left-hand side of \eqref{theoeq} is majorized by 
\begin{equation} \label{Cauchys}
\ll \frac{z}{\Delta^2} \int_{z}^{2z} \sum_{1\leq k\le K} \left|
 \sum_{t< n^2+k \leq t+\Delta} \Lambda(n^2+k) -
 \mathfrak{S}(k) \sum_{t< n^2+k \leq t+\Delta} 1 \right|^2 \dif t.
\end{equation}
Now \eqref{theoeq} follows from Theorem~\ref{sharperresult} and
 \eqref{Cauchys}.
\end{proof}

{\bf Acknowledgements.}
This paper was written when the second-named author was supported by a postdoctoral fellowship at {\it Institutionen f\"or Matematik}, {\it Kungliga Tekniska
H\"ogskolan} in Stockholm and a grant from the G\"oran Gustafsson Foundation.  He wishes to thank these sources for their support.  Furthermore, for the discussions with and encouragements given to him during this work, he owes a debt of gratitude to Prof. John B. Friedlander.

\bibliography{biblio}
\bibliographystyle{amsxport}

\vspace*{.7cm}

\noindent School of Engineering and Science, Jacobs University Bremen \newline
P. O. Box 750561, Bremen 28725 Germany \newline
Email: {\tt s.baier@iu-bremen.de} \newline

\noindent Department of Mathematics, Royal Institute of Technology(KTH) \newline
Lindstedtsv\"agen 25, Stockholm 10044 Sweden \newline
Email: {\tt lzhao@math.kth.se}
\end{document}